\documentclass[12pt]{amsart}
\def\Image{\operatorname{Im}}
\def\Int{\operatorname{Int}}

\newtheorem{Theorem}{Theorem}[section]
\newtheorem{Corollary}[Theorem]{Corollary}
\newtheorem{Lemma}[Theorem]{Lemma}

\theoremstyle{definition}

\newtheorem*{Problem}{Problem}
\theoremstyle{remark}

\begin{document}
\sloppy
\title{Dense subsets of boundaries of CAT(0) groups}
\author{Tetsuya Hosaka} 
\address{Department of Mathematics, Utsunomiya University, 
Utsunomiya, 321-8505, Japan}
\date{September 19, 2005}
\email{hosaka@cc.utsunomiya-u.ac.jp}
\thanks{
Partly supported by the Grant-in-Aid for Young Scientists (B), 
The Ministry of Education, Culture, Sports, Science and Technology, Japan.
(No.~15740029).}
\maketitle
\begin{abstract}
In this paper, 
we study dense subsets of boundaries of CAT(0) groups.
Suppose that a group $G$ acts geometrically on a CAT(0) space $X$ and 
suppose that there exists an element $g_0\in G$ such that 
(1) $Z_{g_0}$ is finite, 
(2) $X\setminus F_{g_0}$ is not connnected, and 
(3) each component of $X\setminus F_{g_0}$ is convex and not $g_0$-invariant,
where $Z_{g_0}$ is the centralizer of $g_0$ and 
$F_{g_0}$ is the fixed-point set of $g_0$ in $X$ 
(that is, $Z_{g_0}=\{h\in G|\, g_0h=hg_0\}$ and $F_{g_0}=\{x\in X|\, g_0x=x\}$).
Then we show that each orbit $G \alpha$ is dense in the boundary $\partial X$ 
(i.e.\ $\partial X$ is minimal) 
and the set $\{g^{\infty}\,|\,g\in G, o(g)=\infty\}$ is also 
dense in the boundary $\partial X$.
We obtain an application for dense subsets on the boundary of a Coxeter system.
\end{abstract}

%%%%%%%%%%%%%
% Section 1 %
%%%%%%%%%%%%%
\section{Introduction and preliminaries}

In this paper, 
we study dense subsets of boundaries of CAT(0) groups.
Definitions and basic properties of CAT(0) spaces and their boundaries 
are found in \cite{BH}.
A {\it geometric} action on a CAT(0) space 
is an action by isometries which is proper (\cite[p.131]{BH}) 
and cocompact.
We note that every CAT(0) space 
on which some group acts 
geometrically is a proper space (\cite[p.132]{BH}).

Suppose that 
a group $G$ acts on a compact metric space $Y$ by homeomorphisms.
Then $Y$ is said to be {\it minimal}, 
if every orbit $Gy$ is dense in $Y$.

For a negatively curved group $G$ and 
the boundary $\partial G$ of $G$, by an easy argument, 
we can show that $G\alpha$ is dense in $\partial G$ 
for each $\alpha\in\partial G$, 
that is, $\partial G$ is minimal.

Suppose that a group $G$ acts geometrically on a CAT(0) space $X$.
Then $G$ acts on the boundary $\partial X$ by homeomorphisms.
The boundary $\partial X$ need not be minimal in general.
Indeed, there exist some examples of Coxeter systems 
whose boundaris are not minimal (cf.\ \cite{Ho1} and \cite{Ho3}).
In \cite{Ho1} and \cite{Ho3}, 
we have investigated dense subsets (and dense orbits) 
of the boundary of a Coxeter system.
In this paper, we study dense subsets of boundaries of CAT(0) groups.
The purpose of this paper is to prove the following theorem.

\begin{Theorem}\label{Thm}
Suppose that a group $G$ acts geometrically on a CAT(0) space $X$.
If there exists an element $g_0\in G$ such that 
\begin{enumerate}
\item[(1)] $Z_{g_0}$ is finite, 
\item[(2)] $X\setminus F_{g_0}$ is not connnected, and 
\item[(3)] each component of $X\setminus F_{g_0}$ is convex and not $g_0$-invariant, 
\end{enumerate}
then each orbit $G \alpha$ is dense in the boundary $\partial X$ 
(i.e.\ $\partial X$ is minimal) 
and the set $\{g^{\infty}\,|\,g\in G, o(g)=\infty\}$ 
is also dense in $\partial X$.
Here $Z_{g_0}$ is the centralizer of $g_0$ and 
$F_{g_0}$ is the fixed-point set of $g_0$ in $X$.
\end{Theorem}

We give the proof of this theorem in Section~2, 
and 
we give some applictions of this theorem and introduce some open problems 
about dense orbits of boundaries of CAT(0) groups 
in Section~3.

%%%%%%%%%%%%%
% Section 2 %
%%%%%%%%%%%%%
\section{Proof of the main theorem}

We first show some results about 
CAT(0) groups and their boundaries needed in the proof of Theorem~\ref{Thm}.

Suppose that a group $G$ acts geometrically on a CAT(0) space $X$.
For an element $g\in G$, 
we define $Z_g$ as the centralizer of $g$ and 
define $F_{g}$ and $\mathcal{F}_{g}$ 
as the fixed-point sets of $g$ in $X$ and $\partial X$ respectively, 
that is, 
\begin{align*}
Z_g&=\{h\in G|\, gh=hg\}, \\
F_{g}&=\{x\in X|\, gx=x\} \ \text{and} \\
\mathcal{F}_{g}&=\{\alpha\in\partial X|\, g\alpha=\alpha\}.
\end{align*}
Also for a subset $A\subset G$, 
the limit set $L(A)$ of $A$ is defined as 
$$L(A)=\overline{Ax_0}\cap \partial X, $$
where $x_0\in X$ and 
$\overline{Ax_0}$ is the closure 
of the orbit $Ax_0$ in $X\cup \partial X$.
We note that 
the limit set $L(A)$ is determined by $A$ and not depend on the point $x_0$.

The following theorem has been proved in \cite{Ho4}.

\begin{Theorem}\label{thm2}
Suppose that a group $G$ acts geometrically on a CAT(0) space $X$.
For $g\in G$, 
the equality $\mathcal{F}_g=L(Z_g)$ holds.
Moreover, 
$\mathcal{F}_g=\emptyset$ if and only if $Z_g$ is finite.
\end{Theorem}

Here 
if $F_g$ is unbounded then $\mathcal{F}_g\neq\emptyset$.
Hence if $Z_g$ is finite then $F_g$ is bounded.

We show the following basic lemma.

\begin{Lemma}\label{lem2}
Suppose that a group $G$ acts geometrically on a CAT(0) space $X$.
Let $A$ be a subset of the boundary $\partial X$, 
let $x_0\in X$ and let $M,N>0$.
If for each $g\in G$ such that $d(x_0,gx_0)>M$, 
there exists $\alpha\in A$ such that $d(gx_0,\Image \xi_{\alpha})\le N$, 
then $A$ is dense in $\partial X$.
Here $\xi_{\alpha}$ is the geodesic ray in $X$ such that 
$\xi_{\alpha}(0)=x_0$ and $\xi_{\alpha}(\infty)=\alpha$.
\end{Lemma}

\begin{proof}
Suppose that 
for each $g\in G$ such that $d(x_0,gx_0)>M$, 
there exists $\alpha\in A$ such that $d(gx_0,\Image \xi_{\alpha})\le N$.

Let $\beta\in\partial X$ and 
let $\xi_{\beta}$ be the geodesic ray in $X$ such that 
$\xi_{\beta}(0)=x_0$ and $\xi_{\beta}(\infty)=\beta$.
To prove that $A$ is dense in $\partial X$, 
we show that 
for each $\epsilon>0$ and $R>0$, 
there exists $\alpha\in A$ such that 
$d(\xi_{\beta}(R),\Image\xi_{\alpha})<\epsilon$.

Let $\epsilon>0$ and $R>0$.
Since the action of $G$ on $X$ is cocompact, 
$B(Gx_0,N')=X$ for some $N'>0$.
Let $R'$ be a number such that 
$$ R'>\max\left\{\frac{R(N+N')}{\epsilon}, M+N'\right\}.$$
Since $\xi_{\beta}(R')\in X=B(Gx_0,N')$, 
there exists $g\in G$ such that $d(\xi_{\beta}(R'),gx_0)\le N'$.
Then 
\begin{align*}
d(x_0,gx_0)&\ge d(x_0,\xi_{\beta}(R'))-d(\xi_{\beta}(R'),gx_0) \\
&\ge R'-N'>(M+N')-N'=M.
\end{align*}
Hence $d(x_0,gx_0)>M$, and 
there exists $\alpha\in A$ such that $d(gx_0,\Image \xi_{\alpha})\le N$.
Then 
$$d(\xi_{\beta}(R'),\Image\xi_{\alpha})\le
d(\xi_{\beta}(R'),gx_0)+d(gx_0,\Image\xi_{\alpha})\le N'+N.$$
Hence
\begin{align*}
d(\xi_{\beta}(R),\Image\xi_{\alpha})
&\le \frac{R}{R'}d(\xi_{\beta}(R'),\Image\xi_{\alpha}) \\
&\le \frac{R}{R'}(N+N')<\epsilon.
\end{align*}
Thus $A$ is dense in $\partial X$.
\end{proof}

The following lemma is known.

\begin{Lemma}[{cf.\ \cite[Proposition~II.6.2(2)]{BH}}]\label{lem3}
Let $X$ be a CAT(0) space and 
let $g$ and $h$ be isometries of $X$.
Then $gF_h=F_{ghg^{-1}}$, 
where $F_h$ is the fixed-point set of $h$ in $X$.
\end{Lemma}

Using the above results, we prove Theorem~\ref{Thm}.

\begin{proof}[Proof of Theorem~\ref{Thm}]
Suppose that a group $G$ acts geometrically on a CAT(0) space $X$ 
and suppose that there exists an element $g_0\in G$ such that 
\begin{enumerate}
\item[(1)] $Z_{g_0}$ is finite, 
\item[(2)] $X\setminus F_{g_0}$ is not connnected, and 
\item[(3)] each component of $X\setminus F_{g_0}$ is convex 
and not $g_0$-invariant.
\end{enumerate}
Then $F_{g_0}$ is bounded by (1) and Theorem~\ref{thm2}.
Let $x_0$ be a basepoint of $X$ and 
let $N$ be a number such that $F_{g_0}\subset B(x_0,N)$.

To prove that every orbit $G \alpha$ is dense in $\partial X$, 
by Lemma~\ref{lem2},
we show that 
for each $\alpha\in\partial X$ and $g\in G$, 
there exists $h\in G$ such that $d(gx_0,\Image \xi_{h\alpha})\le N$, 
where $\xi_{h\alpha}$ is the geodesic ray in $X$ 
such that $\xi_{h\alpha}(0)=x_0$ and $\xi_{h\alpha}(\infty)=h\alpha$.

Let $\alpha\in\partial X$ and $g\in G$.

If $\Image \xi_{\alpha}\cap gF_{g_0}\neq\emptyset$, then 
$$ d(gx_0,\Image\xi_{\alpha})\le d(gx_0,\Image\xi_{\alpha}\cap gF_{g_0})
\le N,$$
because $gF_{g_0}\subset B(gx_0,N)$.

Suppose that $\Image \xi_{\alpha}\cap gF_{g_0}=\emptyset$.
Let $gX'$ be the component of $X\setminus gF_{g_0}$ 
such that $x_0\in gX'$.
Since $\Image \xi_{\alpha}\cap gF_{g_0}=\emptyset$, 
we have that $\Image \xi_{\alpha}\subset gX'$.
Then $gg_0g^{-1}(\Image\xi_{\alpha})\subset gg_0g^{-1}(gX')=gg_0X'$.
Here $gX'\neq gg_0X'$ because $X'\neq g_0X'$ 
by non-$g_0$-invariantness of $X'$ 
which is a component of $X\setminus F_{g_0}$.
Hence $gg_0g^{-1}\alpha\in\partial(gg_0X')$.
On the other hand, $x_0\in gX'$.
Thus $\Image \xi_{gg_0g^{-1}\alpha}\cap gF_{g_0}\neq\emptyset$, 
where $\xi_{gg_0g^{-1}\alpha}$ is the geodesic ray in $X$ 
such that $\xi_{gg_0g^{-1}\alpha}(0)=x_0$ 
and $\xi_{gg_0g^{-1}\alpha}(\infty)=gg_0g^{-1}\alpha$.
Hence 
$$ d(gx_0,\Image\xi_{gg_0g^{-1}\alpha})
\le d(gx_0,\Image\xi_{gg_0g^{-1}\alpha}\cap gF_{g_0})
\le N,$$
since $gF_{g_0}\subset B(gx_0,N)$.

Therefore every orbit $G \alpha$ is dense in $\partial X$.

Next we prove that 
the set $\{h^{\infty}\,|\,h\in G, o(h)=\infty\}$ is dense in $\partial X$.
By Lemma~\ref{lem2}, 
we show that 
for each $g\in G$ such that $d(x_0,gx_0)>2N$, 
there exists $h\in G$ such that $o(h)=\infty$ and 
$d(gx_0,\Image\xi_{h^\infty})\le N$.

Let $g\in G$ such that $d(x_0,gx_0)>2N$.

We first suppose that $[x_0,gx_0]\cap gF_{g_0}=\emptyset$, 
where $[x_0,gx_0]$ is the geodesic from $x_0$ to $gx_0$ in $X$.
Then 
\begin{align*}
d(gF_{g_0},gg_0gF_{g_0})&=d(F_{g_0},g_0gF_{g_0})\\
&=d(g_0^{-1}F_{g_0},gF_{g_0})\\
&=d(F_{g_0},gF_{g_0})\\
&\ge d(x_0,gx_0)-2N \\
&>0, 
\end{align*}
because $g_0F_{g_0}=F_{g_0}$, $F_{g_0}\subset B(x_0,N)$ and $d(x_0,gx_0)>2N$.
Hence $gF_{g_0}\cap gg_0gF_{g_0}=\emptyset$.
Thus 
$(gg_0)^i[x_0,gx_0]\cap (gg_0)^igF_{g_0}=\emptyset$ 
and $(gg_0)^igF_{g_0}\cap (gg_0)^{i+1}gF_{g_0}=\emptyset$ 
for each $i\in\{0,1,2,\ldots\}$.
We consider 
$$Y=X\setminus\bigcup\{(gg_0)^igF_{g_0}\,|\,i=0,1,2,\ldots\}.$$
Let $Y_i$ be the component of $Y$ such that $(gg_0)^ix_0\in Y_i$ 
for each $i\in\{0,1,2,\ldots\}$.
Then $Y_i\neq Y_j$ if $i\neq j$ by the above argument.
Hence $o(gg_0)=\infty$.
We consider $(gg_0)^\infty$ and the geodesic ray $\xi_{(gg_0)^\infty}$.
Since the sequence $\{(gg_0)^i\}$ converges to $(gg_0)^\infty$ and 
$(gg_0)^ix_0\in Y_i$ for each $i$, 
we obtain that $\Image \xi_{(gg_0)^\infty}\cap (gg_0)^igF_{g_0}\neq\emptyset$
for each $i$.
In particular, 
$\Image \xi_{(gg_0)^\infty}\cap gF_{g_0}\neq\emptyset$.
Thus 
$$ d(gx_0,\Image\xi_{(gg_0)^\infty})
\le d(gx_0,\Image \xi_{(gg_0)^\infty}\cap gF_{g_0})
\le N,$$
since $gF_{g_0}\subset B(gx_0,N)$.

Next we suppose that $[x_0,gx_0]\cap gF_{g_0}\neq\emptyset$.
Let $gX'$ be the component of $X\setminus gF_{g_0}$ 
such that $x_0\in gX'$.
Since $X'$ is not $g_0$-invariant, 
there exists a number $k$ such that $gg_0^kx_0\in gX'$.
Then $[x_0,gg_0^kx_0]\cap gF_{g_0}=\emptyset$ 
because $x_0,gg_0^kx_0\in gX'$ and $gX'$ is convex.
Here 
\begin{align*}
d(gg_0^kF_{g_0},(gg_0^k)g_0(gg_0^k)F_{g_0})&=d(F_{g_0},g_0(gg_0^k)F_{g_0})\\
&=d(F_{g_0},g_0gF_{g_0}) \\
&=d(g_0^{-1}F_{g_0},gF_{g_0}) \\
&=d(F_{g_0},gF_{g_0}) \\
&\ge d(x_0,gx_0)-2N \\
&>0.
\end{align*}
Hence 
$gg_0^kF_{g_0}\cap (gg_0^k)g_0(gg_0^k)F_{g_0}=\emptyset$.
By the same argument as the above one, 
we have that 
$\Image\xi_{((gg_0^k)g_0)^\infty}\cap gg_0^kF_{g_0}\neq\emptyset$.
Since $g_0^kF_{g_0}=F_{g_0}$, 
$\Image\xi_{(gg_0^{k+1})^\infty}\cap gF_{g_0}\neq\emptyset$.
Thus 
$$ d(gx_0,\Image\xi_{(gg_0^{k+1})^\infty})
\le d(gx_0,\Image \xi_{(gg_0^{k+1})^\infty}\cap gF_{g_0})
\le N.$$

Therefore 
the set $\{h^{\infty}\,|\,h\in G, o(h)=\infty\}$ is dense in $\partial X$.
\end{proof}

%%%%%%%%%%%%%
% Section 3 %
%%%%%%%%%%%%%
\section{Applications and problems}

In \cite{Ho2}, 
we have defined a {\it reflection} of a geodesic space.
We say that an isometry $r$ of a geodesic space $X$ is a 
{\it reflection} of $X$, if 
\begin{enumerate}
\item[(1)] $r^2$ is the identity of $X$, 
\item[(2)] $X\setminus F_r$ has exactly 
two convex connected components $X^+_r$ and $X^-_r$, 
\item[(3)] $rX^+_r=X^-_r$ and 
\item[(4)] $\Int F_r=\emptyset$, 
\end{enumerate}
where $F_r$ is the fixed-point set of $r$.

We obtain the following corollary from Theorem~\ref{Thm}.

\begin{Corollary}\label{Cor1}
Suppose that a group $G$ acts geometrically on a CAT(0) space $X$.
If there exists a reflection $r\in G$ of $X$ such that $Z_r$ is finite, 
then each orbit $G \alpha$ is dense in the boundary $\partial X$ 
(i.e.\ $\partial X$ is minimal) 
and the set $\{g^{\infty}\,|\,g\in G, o(g)=\infty\}$ 
is dense in $\partial X$.
Here $Z_r$ is the centralizer of $r$.
\end{Corollary}

For a Coxeter system $(W,S)$ and the Davis complex $\Sigma(W,S)$ 
which is a CAT(0) space (\cite{D1} and \cite{M}), 
each $s\in S$ is a reflection of $\Sigma(W,S)$.
Hence we obtain the following corollary.

\begin{Corollary}\label{Cor2}
Let $(W,S)$ be a Coxeter system and 
let $\Sigma(W,S)$ be the Davis complex of $(W,S)$.
If there exists $s\in S$ such that $Z_s$ is finite, 
then each orbit $W \alpha$ is dense in the boundary $\partial\Sigma(W,S)$ 
(i.e.\ $\partial\Sigma(W,S)$ is minimal) 
and the set $\{w^{\infty}\,|\,w\in W, o(w)=\infty\}$ 
is dense in $\partial\Sigma(W,S)$.
\end{Corollary}

The following problems are open.

\begin{Problem}
Does there exist a Coxeter system $(W,S)$ 
such that some orbit $W \alpha$ is dense in $\partial\Sigma(W,S)$
and $\partial\Sigma(W,S)$ is not minimal?
\end{Problem}

\begin{Problem}
Suppose that a group $G$ acts geometrically on two CAT(0) spaces $X$ and $X'$.
Is it the case that $\partial X$ is minimal 
if and only if $\partial X'$ is minimal?
\end{Problem}

\begin{Problem}
Suppose that a group $G$ acts geometrically on a CAT(0) space $X$.
Is it always the case that 
the set $\{g^{\infty}\,|\,g\in G, o(g)=\infty\}$ 
is dense in $\partial X$?
\end{Problem}

%%%%%%%%%%%%%%%%%%%%%%%%%%%%%%%%%%%%%
%             REFERENCES            %
%%%%%%%%%%%%%%%%%%%%%%%%%%%%%%%%%%%%%
%

%
\end{document}